\documentclass[a4paper,12pt]{amsart}
\usepackage[utf8]{inputenc}
\usepackage{amsfonts,amssymb,amsthm,amsmath}
\usepackage{color,url,hyperref}
\theoremstyle{plain}
\newtheorem{theorem}{Theorem}[section]

\newtheorem{lemma}[theorem]{Lemma}
\newtheorem{corollary}[theorem]{Corollary}

\newtheorem{conjecture}[theorem]{Conjecture}

\newtheorem*{problemintro*}{Problem}

\theoremstyle{definition}

\newtheorem{problem}[theorem]{Problem}

\definecolor{darkblue}{rgb}{0,0,0.7} 

\newcommand{\polylog}{\rm{polylog}}

\title{Embedding Divisor and Semi-Prime Testability in $f$-vectors of polytopes}

\author{Eran Nevo}
\thanks{Partially supported by the Israel Science Foundation grant ISF-2480/20 and by ISF-BSF joint grant 2016288}

\begin{document}
\maketitle
\begin{abstract}
We obtain computational hardness results for $f$-vectors of polytopes by exhibiting reductions of the problems DIVISOR and SEMI-PRIME TESTABILITY to problems on $f$-vectors of polytopes. Further, we show that the corresponding problems for $f$-vectors of simplicial polytopes are polytime solvable.
The regime where we prove this computational difference (conditioned on standard conjectures on the density of primes and on $P\neq NP$) is when the dimension $d$ tends to infinity and the number of facets is linear in $d$.

\end{abstract}

\section{Introduction}
The $f$-vector $(f_0(P), f_1(P),\ldots,f_{d-1}(P))$ of a $d$-polytope $P$ records the number of faces $P$ has: $f_i(P)$ faces in dimension $i$.
The $f$-vectors of polytopes of dimension at most $3$ were characterized by Steinitz, and the conditions, which are linear equalities and inequalities on the entries of the $f$-vector, are then easy to check; see e.g.~\cite[Sec.10.3]{Grunbaum:ConvexPolytopes-03}. In contrast, the $f$-vectors of $d$-polytopes for $d\geq 4$ are not well understood; see e.g.~\cite[Sec.10.4]{Grunbaum:ConvexPolytopes-03} and the fatness parameter~\cite{Ziegler-ICM} for $d=4$, while the case $d>4$ is even less understood. The set of $f$-vectors of the important subfamily of \emph{simplicial} polytopes is characterized by the $g$-theorem, conjectured by McMullen~\cite{McMullen:NumberFaces-71} and proved by Stanley~\cite{Stanley:NumberFacesSimplicialPolytope-80} and Billera-Lee~\cite{BilleraLee:SufficiencyMcMullensConditions-1981}.
While this well-understood set may be regarded as complicated from some viewpoints (e.g. it is not a semi-algebraic set of lattice points, for any $d\geq 6$, see~\cite{Sjoberg-Ziegler:semi.alg}), yet deciding membership in it is computationally easy, see~\cite[Thm.1.4]{Nevo-yardsticks}.
The analogous computational problem for the set of $f$-vectors of $d$-polytopes is unsolved, see~\cite[Problem 1.5]{Nevo-yardsticks}, and we conjecture it to be NP-hard. It is known to be decidable in time double exponential in the input size.

We exhibit two variants of the above membership problem and show that they are computationally hard for $f$-vectors of polytopes (given standard conjectures in complexity theory), but are efficiently solvable for $f$-vectors of simplicial polytopes.

\begin{problem}(Fiber Count)\label{prob:fiber-count}
Given $d$, a subset of integers $S\subseteq [0,d-1]$, and values $f_i$ for all $i\in S$, let $\rm{fc}=\rm{fc} (d,(f_i)_{i\in S})$ be the number of $f$-vectors of $d$-polytopes with the given values for the $S$-coordinates. What is the computational complexity:

 (i) of computing $\rm{fc}$ as a function of the input size $N$?

(ii) of deciding whether $\rm{fc}=1$?
\end{problem}

The problem of computing the number of divisors of a given integer, or even of deciding if a given integer is the product of exactly two primes (Semiprime Testability), is believed to be as hard as FACTORING, namely, as factoring the integer into a product of primes; see e.g. Terry Tao's answer at MathOverflow~\cite{Tao-MO}. From a structural result of McMullen on $d$-polytopes with $d+2$ facets~\cite{McMullen:LBT-general_polytopes}, specialized to the case $f_0=2d+1$ (see also~\cite{Pineda:LBT-general_polytopes}), we conclude:
\begin{lemma}\label{lem:LBT-Semiprime}
(i) The number of $f$-vectors of $d$-polytopes with
$f_0=2d+1$ and $f_{d-1}=d+2$ equals $\lceil\frac{D(d)}{2}\rceil$, where $D(d)$ is the number of divisors of $d$ in the interval $[2,d-1]$.

(ii) In particular, $\rm{fc} (d,f_0=2d+1,f_{d-1}=d+2)=1$ iff $d$ is a either a semiprime or equals $p^3$ for some prime $p$.
\end{lemma}

As a corollary, we can reduce Semiprime Testability to a decision problem on fiber count, namely Problem~\ref{prob:fiber-count}(ii).
Here the bit length of the input is $O(\log d)$, while the full $f$-vector clearly has bit length of size $\Omega(d)$ (in fact $\Omega(d \log d)$). Nevertheless, the corresponding problem for $f$-vectors of simplicial polytopes can be solved efficiently:

Let $\rm{fc}_s=\rm{fc}_s (d,(f_i)_{i\in S})$ be the number of $f$-vectors of \emph{simplicial} $d$-polytopes with the given values for the $S$-coordinates.

\begin{theorem}\label{thm:simplicial-FiberCount}
Given as input positive integers $d, a, b$ of total bit length $O(\log d)$,
and $b$ of order $O(d)$:

(i)
It can be decided in $\polylog(d)$-time whether
$\rm{fc}_s (d,f_0=a,f_{d-1}=b) =1$.

(ii)
Deciding whether
$\rm{fc} (d,f_0=a,f_{d-1}=b) =1$ is at least as hard as Semiprime Testability for $d$.
\end{theorem}

The problem DIVISOR, asking whether given three integers $L<U<d$, $d$ has a divisor in the interval $[L,U]$, is believed to be NP-complete,
see e.g. Sudan's survey~\cite{Sudan-PvsNP}. In fact, it is NP-complete if for any large enough real number $x$ there exists a prime in the interval $[x,x+\rm{polylog}(x)]$, see e.g. Peter Shor and Boaz Barak answers at StackExchange~\cite{Shor-StackExchange} to a question by Micha\"{e}l Cadilhac. Cram\'{e}r conjecture~\cite{Cramer36:gap_between_primes, Granville:Cramer} implies that for any $\epsilon >0$ the interval $[x,x+(1+\epsilon)\log^2(x)]$ suffices for $x$ large enough. DIVISOR remains NP-complete if we require $\sqrt{d}\in [L,U]$ (under the assumption above on the existence of primes in short intervals), by a reduction from a variant of SUBSET SUM of real numbers where the target sum is approximately \emph{half} the sum of all input numbers.

\begin{lemma}\label{lem:LBT-DIVISOR}
Given three integers $L<U<d$, with $\sqrt{d}\in [L,U]$, denote $M=M(L,U,d)=\max(L+\frac{d}{L}, U+\frac{d}{U})$. Then there exists a divisor $x$ of $d$ such that $L\le x\le U$ iff there exists a $d$-polytope $P$ whose $f$-vector satisfies $f_0(P)=2d+1$, $f_{d-1}(P)=d+2$ and $f_1(P)\in
[d^2+\frac{d}{2}(1+d-M), d^2+\frac{d}{2}(1+d-2\sqrt{d})]$.
\end{lemma}

Again, we show that the corresponding problem for simplicial polytopes is polytime-solvable, despite the fact that the input is of size logarithmic in $d$, the number of coordinates in the $f$-vector. Combined, it read as follows.

\begin{theorem}\label{thm:simplicial-Range}
Given as input positive integers $d, a, b, L, U$ of total bit length $O(\log d)$, such that $L\leq U$ and
$b$ is of order $O(d)$,
then:

(i)
It can be decided in $\polylog (d)$-time whether there exists a simplicial $d$-polytope $P$ whose $f$-vector satisfies $f_0(P)=a$, $f_{d-1}(P)=b$ and $f_1(P)\in [L,U]$.

(ii)
Deciding whether there exists a $d$-polytope $P$ whose $f$-vector satisfies $f_0(P)=a$, $f_{d-1}(P)=b$ and $f_1(P)\in [L,U]$ is at least as hard as DIVISOR for $d$.
\end{theorem}

Let us remark that Sj\"{o}berg and Ziegler
characterized the pairs $(n,m)$ such that  there exists a $d$-polytope $P$ with $(f_0(P),f_{d-1}(P))=(n,m)$ for even $d$ in the regime
$n+m \ge \binom{3d+1}{\lfloor d/2 \rfloor}$
(and they proved similar but weaker results for $d$ odd); however our interest is in the regime $m+n\in O(d)$ where the behaviour is different and not well understood.

\textbf{Outline.} Section~\ref{sec:Prelim} sets notation and collects the background results we need on $f$-vectors of polytopes.
In Section~\ref{sec:Reductions} we prove the computational hardness results above, for general polytopes, namely Theorems~\ref{thm:simplicial-FiberCount}(ii) and~\ref{thm:simplicial-Range}(ii).
In Section~\ref{sec:eff-simp} we prove the computational efficiency results above, for simplicial polytopes, namely, Theorems~\ref{thm:simplicial-FiberCount}(i) and~\ref{thm:simplicial-Range}(i).
Section~\ref{sec:ConcludingRemarks} ends with open problems.

\section{Preliminaries}\label{sec:Prelim}
For the basics on face enumeration and on polytopes needed here we refer to e.g. the textbooks by Gr\"{u}nbaum~\cite{Grunbaum:ConvexPolytopes-03} and Ziegler~\cite{Ziegler}.

\subsection{Faces of polytopes}
A \emph{$d$-polytope} is a polytope of dimension $d$. Its faces of dimension $k$ are called \emph{$k$-faces}.
Faces of dimension $0$, $1$, $d-1$ are called \emph{vertices, edges, facets}, respectively.
A polytope is \emph{simplicial} if all its proper faces are simplices.

Denote by $f_k(P)$ the number of $k$-faces of a $d$-polytope $P$. The $f$-vector of $P$ is
$f(P)=(1=f_{-1}(P),f_0(P),f_1(P),\ldots,f_{d-1}(P))$.

The following lower bound result of McMullen is crucial for our computational hardness results: let
\[\Phi_j(v,d)=\min\{f_j(P):\ P\ \text{is a d-polytope},\ f_0(P)=v\}
\]

\begin{theorem}\cite[Thm.2]{McMullen:LBT-general_polytopes}\label{thm:McMullen-LBT}
\begin{enumerate}
  \item $\Phi_{d-1}(d+1,d)=d+1$, achieved by the $d$-simplex only.
  \item If $d+2\leq v\leq \lfloor \frac{d(d+8)}{4}\rfloor$, then either (i) $\Phi_{d-1}(v,d)=d+2$, and a $d$-polytope that achieves this must be of the form $T^{r,s,t}:=$ a $t$-fold pyramid over the cartesian product of an $r$-simplex and an $s$-simplex. Thus $v=(r+1)(s+1)+t$, $d=r+s+t$, $t\geq 0$ and $r,s\geq 1$ for some integers $r,s,t$ in this case. (ii) Or else, $\Phi_{d-1}(v,d)=d+3$.
\end{enumerate}
\end{theorem}

\subsection{Face numbers of simplicial polytopes}
Assume that the $d$-polytope $P$ is \emph{simplicial}.
Then the $f$-vector and $h$-vector of $P$ determine each other by a polynomial equation in the ring $\mathbb{Z}[x]$:
\[\sum_{i=0}^{d}f_{i-1}x^{d-i} = \sum_{i=0}^{d}h_{i}(x+1)^{d-i}
.\]
 Define the $g$-vector $g(P)=(g_0,\ldots,g_{\lfloor d/2\rfloor})$ by setting $g_0=1$ and $g_i=h_i-h_{i-1}$ for $1\leq i
\leq d/2$. The celebrated $g$-theorem~\cite{BilleraLee:SufficiencyMcMullensConditions-1981, Stanley:NumberFacesSimplicialPolytope-80} asserts:
\begin{theorem}($g$-theorem)\label{thm:g-thm}
$f=(1,f_0,\ldots,f_{d-1})$ is the $f$-vector of a simplicial $d$-polytope iff

(i) the corresponding $h$-vector satisfies Dehn-Sommerville relations: $h_i=h_{d-i}$ for all $0\leq i\leq \lfloor d/2 \rfloor$; and

(ii) the corresponding $g$-vector is an $M$-sequence, namely $0\leq g_i$ for all $1\leq i\leq d/2$ and it satisfies Macaulay inequalities $g_i^{<i>}\ge g_{i+1}$ for all $1\le i\le \lfloor d/2\rfloor -1$.
\end{theorem}

\section{Reductions}\label{sec:Reductions}
Here we prove our computational hardness results, Theorems~\ref{thm:simplicial-FiberCount}(ii) and~\ref{thm:simplicial-Range}(ii), via Lemmas~\ref{lem:LBT-Semiprime} and~\ref{lem:LBT-DIVISOR} resp.

As observed in~\cite{Pineda:LBT-general_polytopes}, plugging $v=2d+1$ into Theorem~\ref{thm:McMullen-LBT} gives the following, as then $d=sr$.

\begin{corollary}\label{cor:LBT}
\begin{enumerate}
  \item If $d$ is a prime then $\Phi_{d-1}(2d+1,d)=d+3$.
  \item If $d$ is the product of exactly two primes, or equals a prime cubed, then $\Phi_{d-1}(2d+1,d)=d+2$, achieved by a unique minimizer polytope.
  \item If $d$ is the product of more than two primes, and not a prime cubed, then $\Phi_{d-1}(2d+1,d)=d+2$, and is achieved by
      $\lceil \frac{D}{2}\rceil>1$ minimizer polytopes, where $D$ is the number of divisors of $d$ in the interval $[2,d-1]$.  Each of these minimizers have a different number of edges, hence a different $f$-vector.
\end{enumerate}
\end{corollary}
The only part of Corollary~\ref{cor:LBT} that
is not immediate from Theorem~\ref{thm:McMullen-LBT} is the claim on the different $f_1$ in part (3). However, a routine computation gives that $$f_1(T^{r,s,t})=d^2+\frac{d(t+1)}{2}$$
in this case (which is indeed an integer!), hence fixing $f_1$ determines $t$ which in turn determines $r$ and $s$ as $rs=d$ and $r+s=d-t$.

Lemma~\ref{lem:LBT-Semiprime} immediately follows. Theorem~\ref{thm:simplicial-FiberCount}(ii) follows by plugging $a=2d+1$ and $b=d+2$, and recalling that deciding if a given $d$ equals a prime cubed is polytime solvable: first one checks if $d^{1/3}$ is an integer in $O((\log d)^{1+\epsilon})$-time (for any fixed $\epsilon>0$), see e.g.~\cite{D.Bernstein}, and if the answer is Yes, then one checks primality of $d^{1/3}$ in $O(\polylog(d))$-time by~\cite{Prime-is-is-P}.

To prove Lemma~\ref{lem:LBT-DIVISOR} we use again the expression for $f_1(T^{r,s,t})$:
recall we assume that $\sqrt{d}\in [L,U]$.
Note that the function $x\mapsto x+\frac{d}{x}$ has a unique extremal point for $x\geq 0$, which is a local minimum, at $x=\sqrt{d}$. Thus, there exists a divisor $r$ of $d$ with $L\leq r\leq U$ iff there exists $T^{r,s,t}$
with $d-t=r+s=r+\frac{d}{r}\in [2\sqrt{d},M]$ for $M=M(d,L,U):=\max\{L+\frac{d}{L},U+\frac{d}{U}\}$, equivalently with $t\in [d-M,d-2\sqrt{d}]$.
This happens, using Corollary~\ref{cor:LBT}, iff there exists a $d$-polytope $P$ with $f_0(P)=2d+1$, $f_{d-1}(P)=d+2$ and $f_1(P)\in [d^2+\frac{d}{2}(1+d-M), d^2+\frac{d}{2}(1+d-2\sqrt{d})]$; as claimed.

As before, Theorem~\ref{thm:simplicial-Range}(ii) follows from the case $a=2d+1$ and $b=d+2$.

\section{Efficient computations for simplicial polytopes}\label{sec:eff-simp}
Here we prove our computational efficiency  results, Theorems~\ref{thm:simplicial-FiberCount}(i) and~\ref{thm:simplicial-Range}(i) using the $g$-theorem.

By a direct computation, the number of facets is expressed in terms of the $g$-vector as follows: for $d=2k$ even
\[
f_{d-1}=(d+1)+(d-1)g_1+(d-3)g_2+\ldots +3g_{k-1}+g_k
,\]
and for $d=2k+1$ odd
\[
f_{d-1}=(d+1)+(d-1)g_1+(d-3)g_2+\ldots +4g_{k-1}+2g_k
.\]

Now, combined with the $g$-theorem, if $f_{d-1}(P)=b\in O(d)$ then there exists a constant $C>0$ s.t. $g_i(P)=0$ for all $i>C$ and $0\leq g_i(P)\leq C$ for all $0\leq i\leq \lfloor d/2\rfloor$;
hence, there are only finitely many potential $g$-vectors to check. In each of them the Macaulay inequalities $g_i^{<i>}\geq g_{i+1}$ need to be checked only for $i<C$, so each such inequality is checked in constant time.
Altogether, in constant time all the $g$-vectors whose $f_{d-1}$ equals $b$ are found.

In particular, one checks in constant time if there exists exactly one such $g$-vector; this   proves Theorem~\ref{thm:simplicial-FiberCount}(i).

Now, for each $g$-vector which passed the test above we compute $f_1=g_2 + dg_1 + \binom{d+1}{2}$ in $O(\polylog(d))$-time and then check whether $f_1\in [L,U]$ in $O(\log(d))$-time, proving Theorem~\ref{thm:simplicial-Range}(i).

\section{Concluding remarks}\label{sec:ConcludingRemarks}
For fixed dimension we conjecture the following, which may be viewed as an explanation why when $d\geq 4$ the $f$-vectors of $d$-polytopes are poorly understood.
\begin{conjecture}
Let $d\ge 4$ be fixed.
Then it is NP-hard to decide if a given $N$-bit vector $f=(1,f_0,\ldots,f_{d-1})$ of positive integers is the $f$-vector of a $d$-polytope.
\end{conjecture}

Regarding the computational efficiency results,
\begin{problem}
Can the assumption $b\in O(d)$ in
Theorems~\ref{thm:simplicial-FiberCount}(i) and~\ref{thm:simplicial-Range}(i) be dropped and the same conclusions there hold?
\end{problem}
This means $b$ is polynomial (rather than linear) in $d$, as the entire input is of size $O(\log d)$.

\textbf{Acknowledgements.}
I deeply thank Nathan Keller for pointing me to~\cite{D.Bernstein} and~\cite{Tao-MO}, and
to Guillermo Pineda Villavicencio for helpful comments on an earlier version.
\bibliographystyle{plain}
\bibliography{biblioERC}
\end{document}